\documentclass[12pt]{amsart}
\usepackage{a4wide}
\usepackage{amsthm,amsfonts,amsmath,mathrsfs,amssymb}

\usepackage{t1enc}

\newtheorem{theorem}{Theorem}[section]

\newtheorem{lemma}[theorem]{Lemma}

\newcommand{\T}{\mathbb{T}}

\begin{document}
\title{The Kadets 1/4 theorem for polynomials}
\author {Jordi Marzo}
\address{Universitat de Barcelona, Gran Via 585, 08071 Barcelona, Spain} \email{jmarzo@ub.edu}
\author{Kristian Seip}
\address{Department of Mathematical Sciences, Norwegian University
of Science and Technology, NO-7491 Trondheim, Norway}
\email{seip@math.ntnu.no}
\thanks{The first author is supported by projects
2005SGR00611 and MTM2005-08984-C02-02.}
\thanks{The second author is
supported by the Research Council of Norway grant 160192/V30.}
\date{\today}
\begin{abstract}
We determine the maximal angular perturbation of the $(n+1)$th
roots of unity permissible in the Marcinkiewicz--Zygmund theorem
on $L^p$ means of polynomials of degree at most $n$. For $p=2$,
the result is an analogue of the Kadets 1/4 theorem on
perturbation of Riesz bases of holomorphic exponentials.
\end{abstract}

\maketitle

\section{Introduction}                                      \label{SecPrelim}

A classical theorem of J. Mar\-cin\-kiewicz and A. Zygmund
generalizes the elementary mean value formula \begin{equation}
\label{plancherel} \frac{1}{n+1}\sum_{j=0}^n
\left|P\left(e^{i\frac{2\pi j}{n+1}}\right)\right|^2
=\int_0^{2\pi} |P(e^{i\theta})|^2\,
\frac{d\theta}{2\pi},\end{equation} valid for holomorphic
polynomials $P$ of degree at most $n$, in the following way: For
$1<p<\infty$, there is a constant $C_p$ independent of $n$ such
that
\begin{equation} \label{CZ} \frac{C^{-1}_p}{n+1}\sum_{j=0}^n
\left|P\left(e^{i\frac{2\pi j}{n+1}}\right)\right|^p \le
\int_0^{2\pi} |P(e^{i\theta})|^p \, \frac{d\theta}{2\pi}\le
\frac{C_p}{n+1}\sum_{j=0}^n \left|P\left(e^{i \frac{2\pi
j}{n+1}}\right)\right|^p\end{equation} for every complex polynomial
$P$ of degree at most $n$. (See \cite{MZ37} or Theorem 7.5 in
Chapter X of \cite{Zyg68}.) It is natural to ask if the norm
equivalence expressed by \eqref{CZ} remains valid if we replace the
$(n+1)$th roots of unity $\omega_{nj}=e^{i \frac{2\pi j}{n+1}}$ by
$n+1$ points $z_{nj}$ on the unit circle with a less regular
distribution. C. K. Chui, X.-C. Shen, and L. Zhong \cite{CSZ93}
considered this problem and found that the norm equivalence is
stable under small perturbations of the points $\omega_{nj}$. We
will prove the following sharp version of their result:
\begin{theorem} Suppose $1<p<\infty$ and set $q=\max(p,p/(p-1))$.
The following statement holds if and only if $\delta<1/(2q)$: There
is a constant $C_p$ independent of $n$ such that if
$|\arg(z_{nj}\overline{\omega_{nj}})|\le 2\pi \delta/(n+1)$ for
$0\le j\le n$, then \begin{equation} \label{inequality}
\frac{C^{-1}_p}{n+1}\sum_{j=0}^n |P(z_{nj})|^p \le \int_0^{2\pi}
|P(e^{i\theta})|^p\, \frac{d\theta}{2\pi}\le
\frac{C_p}{n+1}\sum_{j=0}^n |P(z_{nj})|^p\end{equation} for every
holomorphic polynomial $P$ of degree at most $n$.
\end{theorem}

We will see that this theorem is a consequence of a general result
of Chui and Zhong \cite{CZ99}, characterizing the so-called $L^p$
Marcinkiewicz--Zygmund families (to be defined below) in terms of
Muckenhoupt $(A_p)$ weights.

Readers familiar with Paley--Wiener spaces will see the analogy with
the Kadets $1/4$ theorem on perturbations of Riesz bases of complex
exponentials in $L^2$ of an interval \cite{Kad64}. One may view
polynomials as discrete versions of band-limited functions, with the
degree of the polynomial being the counterpart to the notion of
``bandwidth''. The identity \eqref{plancherel} is the discrete
analogue of the Plancherel identity or---what amounts to the
same---the Shannon formula for bandlimited functions. In the case
when $p=2$ and we require $\delta<1/4$, our theorem corresponds
precisely to the Kadets $1/4$ theorem. The $L^p$ version
($1<p<\infty$) of the Kadets theorem, analogous to our theorem, can
be found in \cite{LS97}.

It is interesting to note that our problem as well as that of the
classical Kadets theorem fits into a general theory of unconditional
bases in so-called model spaces. (See \cite{Pa79}, \cite{Ni80}, and
\cite{KNP79} for original work and \cite{Ni02} or \cite{Se04} for
more recent expositions.) In particular, the theorem of Chui and
Zhong to be used in this note can be obtained from a theorem given
in \cite{KNP79}. We refer to \cite{M08} for the details of this link
and to \cite{OS05}, where the connection between
Marcinkiewicz--Zygmund inequalities and model spaces was first
mentioned explicitly.

For $p=2$, the proof to be given below is an adaption of
S.~Khrushchev's proof of the classical Kadets $1/4$ theorem
\cite{K79}, and, for general $p$, we act in a similar way as was
done in \cite{LS97}. Khrushchev also showed how to obtain other
perturbation results, such as a theorem of S. Avdonin \cite{A74}.
We will confine ourselves to proving the theorem stated above and
refer to \cite{M08} for the counterpart of Avdonin's theorem as
well as other analogues of results for Paley-Wiener spaces and
families of complex exponentials.

\section{Preliminaries}
Suppose that for each nonnegative integer $n$ we are given a set
$\mathcal{Z}(n)=\{ z_{n j}\}_{j=0}^{n}$ of $n+1$ distinct points on
the unit circle. We denote by $\mathcal{Z}=\{ \mathcal{Z}(n)
\}_{n\ge 0}$ the corresponding triangular family of points. The
family $\mathcal{Z}$ is declared to be uniformly separated if there
exists a positive number $\varepsilon$ such that
    $$\inf_{j\neq k}|z_{nj}-z_{nk}|\ge \frac{\varepsilon}{n+1}$$
    for every $n\ge 0.$

We will say that $\mathcal{Z}$ is an $L^{p}$ Marcinkiewicz--Zygmund
family if there exists a constant $C_p>0$ such that for every $n\ge
0$ and complex polynomial $P$ of degree at most $n$, we have
\begin{equation}    \label{ineq-def-MZ}
    \frac{C_p^{-1}}{n+1}\sum_{j=0}^{n}|P(z_{n j})|^{p}
    \le \int_{0}^{2\pi }|P(e^{i\theta})|^{p}d\theta
    \le \frac{C_p}{n+1}\sum_{j=0}^{n}|P(z_{nj})|^{p}.
\end{equation}
In order to describe such families, we associate with $\mathcal{Z}$
the following generating polynomials
\[ F_n(z)=\prod_{j=0}^n\left(1-\frac{n}{n+1}\overline{z_{jn}}z\right). \]
The theorem of Chui and Zhong reads as follows \cite{CZ99}.
\begin{theorem}                                                     \label{ChuiZhong}
Suppose $1<p<\infty$. The family $\mathcal{Z}=\{
\mathcal{Z}(n)\}_{n\ge 0}$ of points on the unit circle is an
$L^{p}$ Marcinkiewicz--Zygmund family if and only if it is
uniformly separated and there exists a constant $K_p$ such that
\begin{equation} \label{Ap}
    \left( \frac{1}{|I|}\int_{I}|F_{n}(e^{i\theta})|^{p}d\theta\right)^{1/p}
    \left( \frac{1}{|I|}\int_{I}|F_{n}(e^{i\theta})|^{-p/(p-1)}d\theta\right)^{(p-1)/p}\le K_{p}
\end{equation}
    for every subarc $I$ of the unit circle and every $n\ge 0.$
\end{theorem}
In other words, the sequence $|F_n|^p$ satisfies a uniform $(A_p)$
condition.

In the proof of the positive part of the $p=2$ case of our
theorem, we will make use of the equivalence between the $(A_2)$
and Helson--Szeg\H{o} conditions. We will derive the result for
$p\neq 2$ from the $p=2$ case using the following estimate.
\begin{lemma} \label{from2top}
Let $\alpha,\kappa>0$ be given, and set
$\rho_{\kappa n}=\max(1/2,1-\kappa/(n+1))$. If a given triangular family of
real numbers $\delta_{nj}$ satisfies $\sup_{nj}|\delta_{nj}|<1/2$,
then
$$\left| \prod_{j=0}^{n}(z-\rho_{\kappa n}e^{\frac{2\pi i (j+\alpha \delta_{nj})}{n+1}})\right|
= R_n(z) \left| \prod_{j=0}^{n}(z-\rho_{\kappa n}e^{\frac{2\pi i
(j+\delta_{nj})}{n+1}})\right|^{\alpha},$$ where $R_n(z)$ is bounded
from above and below by positive constants, independently of $z\in
\mathbb{T}$ and $n\ge 0$.
\end{lemma}

\begin{proof}
Set
$$P_{\beta}(\theta)=\left| \prod_{j=0}^{n}(e^{i\theta}-\rho_{\kappa n}e^{i \lambda_{j}(\beta)})\right|
,\;\;\mbox{where}\;\; \lambda_{j}(\beta)=\frac{2\pi
j}{n+1}+\frac{2\pi \beta \delta_{nj}}{n+1}.$$  We have
    $$\log P_{\beta}(\theta)-\log P_{0}(\theta)=
    \mbox{Re} \sum_{j=0}^n\int_{\rho_{\kappa n}e^{i\lambda_{j}(0)}}^{\rho_{\kappa n}e^{i\lambda_{j}(\beta)}}\frac{d\xi}{\xi-e^{i\theta}}
    =\sum_{j=0}^n\int_{\lambda_{j}(0)}^{\lambda_{j}(\beta)}h(\theta-t)dt,$$
where
    $$h(t)=\frac{\rho_{\kappa n}\sin t}{1+\rho_{\kappa n}^{2}-2\rho_{\kappa n}\cos t}.$$
By the fundamental theorem of calculus,
    $$\log P_{\beta}(\theta)-\log P_{0}(\theta)=
    \sum_{j=0}^n(\lambda_{j}(\beta)-\lambda_{j}(0))h(\theta-\lambda_{j}(0))+
    \sum_{j=0}^n\int_{\lambda_{j}(0)}^{\lambda_{j}(\beta)}
    \int_{\lambda_{j}(0)}^{t}h'(\theta-\tau)d\tau dt.$$
We compute $h'(t)$ and find that the absolute value of the latter
sum is bounded independently of $\theta$ and $n.$ Therefore,
\begin{align*}
 \log P_{\alpha}(\theta)-\log P_{0}(\theta) &
 =\alpha\sum_{j=0}^n\frac{2\pi\delta_{nj}}{n+1}h(\theta-\lambda_{j}(0))
    +b_{n,\alpha}(z)
    \\
    &
    =\alpha (\log P_{1}(\theta)-\log P_{0}(\theta)-b_{n,1}(z))+b_{n,\alpha}(z)
\end{align*}
with uniform bounds on the $L^\infty$ norms of
$b_{n,\alpha}$. This gives the result because $P_{0}(\theta)$ is
trivially bounded from above and below by positive constants,
independently of $z\in \mathbb{T}$ and $n\ge 0$.
\end{proof}

\section{Proof of the theorem: Sufficiency}

For each set $\mathcal{Z}(n)$, we define $C_{p}(\mathcal{Z}(n))$ as
the minimum of all positive numbers $C$ such that
\[ \frac{C^{-1}}{n+1}\sum_{j=0}^n |P(z_{nj})|^p \le \int_0^{2\pi}
|P(e^{i\theta})|^p\, \frac{d\theta}{2\pi}\le
\frac{C}{n+1}\sum_{j=0}^n |P(z_{nj})|^p\] for every complex
polynomial $P$ of degree at most $n$. Among all sets
$\mathcal{Z}(n)$ satisfying $|\arg(z_{nj}\overline{\omega_{nj}})|\le
2\pi \delta/(n+1)$ for $0\le j\le n$, we may choose a set with
maximal $C_{p}(\mathcal{Z}(n))$. From now on, we will assume that
the points $z_{nj}=\omega_{nj}e^{\frac{2\pi i \delta_{nj}}{n+1}}$ constitute a
set of points with this extremal property. It suffices to show that
the corresponding triangular family is an $L^p$
Marcinkiewicz--Zygmund family. Clearly, this family is uniformly
separated when $\delta<1/(2q)$.

When $p=2$, condition \eqref{Ap} is equivalent to the following
uniform Helson-Szeg\H{o} condition: There exist sequences $u_n$ and
$v_n$ of real functions in $L^\infty(\T)$ such that
\begin{equation}\label{HS} |F_n|^2=e^{u_n+\widetilde{v_n}} \ \
\text{with} \ \  \sup_n \|u_n\|_\infty<\infty \ \ \text{and} \ \
\sup_n \|v_n\|_\infty < \pi/2. \end{equation} Here $v\mapsto
\widetilde{v}$ denotes the conjugation operator.

We need two steps in order to identify the appropriate functions
$u_n$ and $v_n$. In the first step,  we ``pull'' the points $z_{nj}$
more deeply into the unit disc. For $\kappa>0$, we set
$\rho_{\kappa n}=\max(1/2,1-\kappa/(n+1))$. We define
\[ F_{\kappa n}(z)=\prod_{j=0}^n (1-\rho_{\kappa n}\overline{z_{nj}} z). \]
For fixed $\kappa>0$, we find that
\[ |F_{n}(e^{i t})|^2=e^{u_{\kappa n}(e^{it})}|F_{\kappa n}(e^{i t})|^2, \]
with $\sup_n \|u_{\kappa n}\|_\infty <\infty$.

We now move to the second step. Writing
\[ B_{\kappa n}(z)=\prod_{j=0}^n \frac{z-\rho_{\kappa n} z_{nj}}{1-\rho_{\kappa n}\overline{z_{nj}}
z}, \] we get \[ B_{\kappa n}(z)=z^{n+1} \frac{\overline{F_{\kappa
n}(z)}}{F_{\kappa n}(z)}=z^{n+1} \frac{|F_{\kappa
n}(z)|^2}{F^2_{\kappa n}(z)}\] for $|z|=1$. Since $F^2_{\kappa n}$
is an outer function with $F^2_{\kappa n}(0)=1$, this means that
$F^2_{\kappa n}=e^{\widetilde{v_{\kappa n}}}$, where
\[ v_{\kappa n}(e^{i\theta})=\int_{0}^{\theta}\sum_{j=0}^{n}
\frac{1-\rho_{\kappa n}^{2}}{|e^{i\eta}-\rho_{\kappa n}z_{nj}|^{2}}
    d\eta-(n+1)\theta-c\]
    and $c$ is any suitable constant. If we set
\[ c=\sum_{j=0}^{n}\int_{-2\pi\delta_j}^{0}\frac{1-\rho_{\kappa n}^{2}}{|e^{i\eta}-\rho_{\kappa n}\omega_{nj}|^{2}}
    d\eta, \]
then we may write
\[ v_{\kappa n}(e^{i\theta})=\sum_{j=0}^{n}\int_{0}^{\theta-2\pi\delta_j}
\frac{1-\rho_{\kappa n}^{2}}{|e^{i\eta}-\rho_{\kappa n}\omega_{nj}|^{2}}
    d\eta-(n+1)\theta.\]
On the other hand, using that
\[ \int_{\theta}^{\theta+2\pi/(n+1)}\sum_{j=0}^{n}
\frac{1-\rho_{\kappa n}^{2}}{|e^{i\eta}-\rho_{\kappa n}\omega_{nj}|^{2}}d\eta=2\pi
\ \ \text{and} \ \ \left|\sum_{j=0}^{n}
\frac{1-\rho_{\kappa n}^{2}}{|e^{i\eta}-\rho_{\kappa n}\omega_{nj}|^{2}}-(n+1)\right|\le
\frac{C(n+1)}{\kappa},\] we get
\[ v_{\kappa n}(e^{i\theta})=\sum_{j=0}^{n}\int_{\theta}^{\theta-2\pi\delta_j}
\frac{1-\rho_{\kappa n}^{2}}{|e^{i\eta}-\rho_{\kappa n}\omega_{nj}|^{2}}
    d\eta + O(\kappa^{-1})\]
    when $\kappa\to\infty$. Consequently,
\[ \|v_{\kappa n}\|_\infty\le \sup_\theta \int_{\theta}^{\theta+2\pi\delta/(n+1)}
\sum_{j=0}^{n}
\frac{1-\rho_{\kappa n}^{2}}{|e^{i\eta}-\rho_{\kappa n}\omega_{nj}|^{2}}
    d\eta + O(\kappa^{-1})=2 \pi \delta + O(\kappa^{-1}). \]
Assuming $\delta<1/4$, we now obtain \eqref{HS} by choosing
$\kappa$ sufficiently large.

Finally, we consider the case $p\neq 2.$ We introduce the triangular
family given by the sets
    $$\mathcal{Z}_{q/2}(n)=\{e^{i\lambda_{nj}(q/2)}\}_{j=0}^{n}\;\;\mbox{with}\;\;
    \lambda_{nj}(q/2)=\frac{2\pi j}{n+1}+\frac{\pi q \delta_{nj}}{n+1}.$$
If $\delta<1/(2q)$, then the $p=2$ case applies. In other words,
if we set
\[ G_n(z)=\prod_{j=0}^n \left(1-\rho_{\kappa n}
e^{-i\lambda_{nj}(q)}z\right), \] then the functions $|G_n|^2$
meet the uniform $(A_2)$ condition. By Lemma~\ref{from2top} and
H\"{o}lder's inequality, this implies that the functions $|F_n|^p$
satisfy the uniform $(A_p)$ condition.

\section{Proof of the theorem: Necessity}

We will consider the sets
\[\mathcal{Z}(2n)=\left\{ e^{2\pi i j/(2n+1)} \right\}_{j=0}^{n}\bigcup
\left\{ e^{-2\pi i (j-2\delta)/(2n+1)} \right\}_{j=1}^{n},\] which
can be viewed as perturbations of the rotated $(2n+1)$th roots of
unity $e^{2\pi \delta/(2n+1)}\omega_{(2n)j}$. Let $F_{2n}$ be the
generating polynomial for $\mathcal{Z}(2n)$. We set
$\phi_{n}(z)=F_{2n}(z)/(z^{2n+1}-\rho_{2n}^{2n+1})$ and observe that
we may write
\[ \phi_n(z)    =\prod_{j=1}^{n} \frac{
    z-\rho_{2n}e^{\frac{-2\pi i (j-2\delta)}{2n+1}}}{
    z-\rho_{2n}e^{\frac{-2\pi i j}{2n+1}}} .\]
We have
\[
    \log |\phi_{n}(z)|=\Re (\log \phi_{n}(z))
=   \Re \sum_{j=1}^{n} \int_{\Gamma_{nj}}\frac{d\xi}{\xi -z},
\]
where $\Gamma_{nj}$ is the arc with the parametrization
$\Gamma_{nj}(t)=\rho_{2n}e^{\frac{-2\pi i
j}{2n+1}}e^{\frac{it}{2n+1}}$,  $0\le t\le 4\delta \pi.$ It follows
that
$$|\phi_{n}(e^{it})|\longrightarrow \left| \frac{1-e^{it}}{1+e^{it}}\right|^{2\delta}$$
for $0<t<\pi.$ By Fatou's lemma,
\[
    \left(\int_{0}^{\pi}\left|\frac{1-e^{it}}{1+e^{it}}\right|^{2\delta p}dt\right)
    \left(\int_{0}^{\pi}\left|\frac{1-e^{it}}{1+e^{it}}\right|^{-\frac{2\delta p}{p-1}}dt\right)^{p-1}
     \le \liminf_{n} \int_{0}^{\pi}|\phi_{n}|^{p}\left( \int_{0}^{\pi}|\phi_{n}|^{-\frac{p}{p-1}}
     \right)^{p-1}.
\]
Hence, when $\delta=1/2q$, the weights $|\phi_n|^p$ do not meet the
uniform $(A_p)$ condition, and the same holds for the weights
$|F_{2n}|^2$ since the polynomials $z^{2n+1}-\rho_{2n}^{2n+1}$ are
uniformly bounded away from $0$ for $|z|=1$.

\end{document}